\documentclass[12pt, a4paper]{amsart}

\usepackage[english]{babel}
\usepackage{amsmath,enumerate, amsfonts, amssymb,amsthm, xypic}
\usepackage[dvips]{graphics} 

\newtheorem{thm}{Theorem}[section]
\newtheorem{lemma}[thm]{Lemma}
\newtheorem{prop}[thm]{Proposition}
\newtheorem{cor}[thm]{Corollary}

\theoremstyle{definition}
\newtheorem{defi}[thm]{Definition}

\newtheorem{ex}[thm]{Example}
\newtheorem{rmk}[thm]{Remark}

\newcommand{\as}{arc-symmetric }
\newcommand{\AS}{\mathcal{AS}}
\newcommand{\R}{\mathbb R}
\newcommand{\C}{\mathbb C}
\newcommand{\Z}{\mathbb Z}
\newcommand{\f}{\mathbb F_2}
\newcommand{\dans}{\longrightarrow}
\newcommand{\veps}{$G$-equivariant virtual Poincar\'e series }

\DeclareMathOperator{\jac}{jac}
\DeclareMathOperator{\mult}{mult}
\DeclareMathOperator{\Reg}{Reg}

\DeclareMathOperator{\Hom}{Hom}

\DeclareMathOperator{\Bl}{Bl}
\DeclareMathOperator{\im}{im}

\title{Equivariant Virtual Betti numbers}

\author{Goulwen Fichou}
\date{\today}
\address {Institut Math\'ematiques de Rennes, Universit\'e de Rennes 1, Campus de Beaulieu, 35042 Rennes Cedex, France}

\email{goulwen.fichou@univ-rennes1.fr}
\begin{document}

\begin{abstract} We define a generalised Euler characteristic for arc-symmetric sets endowed with a group action. It coincides with the Poincar\'e series in equivariant homology for compact nonsingular sets, but is different in general. We put emphasis on the particular case of $\Z/2\Z$, and give an application to the study of the singularities of Nash function germs via an analog of the motivic zeta function of Denef and Loeser.
\end{abstract}

\maketitle
%%%%%%%%%%%%%%%%%%%%%%%%%%%%%%%%%%%%%%%%%%%%%%%%%%%%%%%%%%%%%%%%%%%%%%%%
%%%%%%%%%%%%%%%%%%%%%%%%%%%%%%%%%%%%%%%%%%%%%%%%%%%%%%%%%%%%%%%%%%%%%%%%%%%%%%%%%
%%%%%%%%%%%%%%%%%%%%%%%%%%%%%%%%%%%%%%%%%%%%%%%%%%%%%%%%%%%%%%%%%%%%%%%%%%%%%%%%%
%%%%%%%%%%%%%%%%%%%%%%%%%%%%%%%%%%%%%%%%%%%%%%%%%%%%%%%%%%%%%%%%%%%%%%%%%%%%%%%%%

\section{Introduction}
                                
Let $G$ be a finite group. The Grothendieck group $K^0(\mathcal V_{\R}(\R),G)$ of real algebraic varieties with a $G$-action is the abelian group generated by symbols $[X,G]$ where $X$ represents a $G$-equivariant isomorphism class of real algebraic varieties endowed with an algebraic $G$-action, subject to the relation $[X,G]=[Y,G]+[X\setminus Y,G]$ if $Y$ is a closed subvariety of $X$ stable under the $G$-action on $X$. We construct a group homomorphism $\beta^G(\cdot):K^0(\mathcal V_{\R}(\R),G) \dans \Z[[u^{-1}]][u]$, the virtual $G$-equivariant Poincar\'e series (theorem \ref{main}). For $X$ compact and nonsingular, the coefficients of $\beta^G(X)$ coincide with the $G$-equivariant Betti numbers with coefficients in $\f$ (cf. \cite{VH} for the definition of equivariant homology, or section \ref{sect-equiv} below). This \veps is the first non-trivial homomorphism defined on $K^0(\mathcal V_{\R}(\R),G)$. The particular case of the trivial group gives back the virtual Poincar\'e polynomial defined by McCrory and Parusi\'nski \cite{MCP}.
Similarly to the non-equivariant case, we define the \veps for the larger Grothendieck group $K^0(\AS,G)$ of \as sets endowed with a $G$-action (cf. \cite{KK1,fichou}). Arc-symmetric sets form a constructible category, i.e. stable under boolean operations, generated after resolution of singularities by connected components of real algebraic sets (cf. proposition \ref{ecl} for a more precise statement).

As an application, we introduce equivariant zeta functions with signs for a Nash function germ, i.e. an analytic germ with semi-algebraic graph, in the spirit of the motivic zeta function of Denef and Loeser \cite{DL1}. The motivic zeta function is defined by considering the motivic measure of some spaces of arcs related to a given function. This measure takes its values in the Grothendieck group of varieties endowed with an action of some group of roots of unity. The \veps allows us to define real analogs for Nash function germs with $G=\{1,-1\}$. It generalises the preceding ``naive'' zeta functions, with or without signs, considered in \cite{KP,fichou,fichou2}. The latter are invariants for the blow-Nash equivalence between Nash function germs. We prove that these new equivariant zeta functions with signs are again invariants for this blow-Nash equivalence. Moreover, we manage to establish a relation between the naive zeta function and the equivariant zeta functions with signs for constant sign germs, solving a problem in \cite{fichou}.

To deal with this issue in section \ref{bN}, we construct in section \ref{sect-eq-virt} the virtual $G$-equivariant Poincar\'e series defined on the Grothendieck group of \as sets endowed with an algebraic $G$-action. This construction is based on some properties of $G$-equivariant homology that we recall in section \ref{sect-equiv}, together with resolution of singularities \cite{hiro} and the weak factorisation theorem \cite{weak}. We prove in particular that, even if the series is not bounded from below, then it has a finite degree equal to the dimension of the set. We focus on the particular cases of the trivial and free actions and give some examples in the interesting case $G\simeq\Z/2\Z$. Then we study in section \ref{G=2} some properties of the virtual $\Z/2\Z$-equivariant Poincar\'e series, relating notably its value to the value of the virtual Poincar\'e polynomial of the quotient (proposition \ref{cle}). We are led in particular to investigate conditions for the existence, as an \as set, of the quotient of an \as set (proposition \ref{eq-lisse}). We specify also the negative part of the series in terms of the fixed point set (proposition \ref{neg}), similarly to the case of $\Z/2\Z$-equivariant homology \cite{VH}.\\

{\bf Acknowledgements.} The author wish to thank C. Mc Crory, J. Van Hamel and S. Maugeais for useful discussions.

%%%%%%%%%%%%%%%%%%%%%%%%%%%%%%%%%%%%%%%%%%%%%%%%%%%%%%%%%%%%%%%%%%%%%%%%%%%%%%%%%
%%%%%%%%%%%%%%%%%%%%%%%%%%%%%%%%%%%%%%%%%%%%%%%%%%%%%%%%%%%%%%%%%%%%%%%%%%%%%%%%%
%%%%%%%%%%%%%%%%%%%%%%%%%%%%%%%%%%%%%%%%%%%%%%%%%%%%%%%%%%%%%%%%%%%%%%%%%%%%%%%%%
\section{Equivariant homology}\label{sect-equiv}
In this section, we recall the definition and some properties of equivariant Borel-Moore homology and equivariant cohomology that will be useful in the sequel (see \cite{VH} for a complete treatment). We do not explain these theories in great generality but mainly focus on the situation we will be interested in in section \ref{sect-eq-virt}.
\subsection{Cohomology of groups}

To begin with, let us recall the definition of the cohomology of a group $G$ with coefficients in an abelian group $M$ endowed with an action of $G$.

Define the integral group ring of $G$ to be $\Z[G]=\{\sum_{g \in G} n_g .g |~n_g \in \Z\}$. An abelian group $M$ with a $G$-action is naturally a $\Z[G]$-module. In order to define the cohomology groups of $G$ with coefficients in $M$, take a projective resolution $(F_{\cdot},\partial_{\cdot})$ of $\Z$ via $\Z[G]$-modules. Then, for $n \in \Z$, define the $n$-th cohomology group of $G$ with coefficients in $M$ by
$$H^n(G,M)=H^n(K^{\cdot},\delta^{\cdot}),$$
where $K^{\cdot}$ is the chain complex defined by $K^i=\Hom _G(F_i,M)$ and $\delta^i(u)=(-1)^i(u \circ \partial _i)$.

\begin{ex} In case $G$ is finite cyclic with generator $\sigma$ and order $d$, one can choose
$$\xymatrix{ \cdots \ar[r] & \Z[G] \ar[r]^{\sigma -1} &  \Z[G] \ar[r]^N  & \Z[G] \ar[r] ^{\sigma -1} & \Z \ar[r]  & 0 }$$
with $N=\sum_{1 \leq i\leq d}\sigma ^i$ as a projective resolution.
Therefore
\begin{displaymath}
H^n(G,M) = \left\{ \begin{array}{lllll}
\frac{\ker (M \xrightarrow{\sigma -1} M)}{\im (M \xrightarrow{N} M)}     & \textrm{if n is even,~~} n>0,\\
\\
\frac{\ker (M \xrightarrow{N} M)}{\im (M \xrightarrow{\sigma -1}M)} & \textrm{if n is odd,~}n> 0,\\
\\
\ker (M \xrightarrow{\sigma -1}M)   & \textrm{if } n = 0.
\end{array} \right.
\end{displaymath}
If $G=\Z/2\Z$, then 
\begin{displaymath}
H^n(G,M) = \left\{ \begin{array}{ll}
 M^G   & \textrm{if~~}n=0 ,\\
\frac{M^G}{\im (1+\sigma)} & \textrm{if~~}n> 0.
\end{array} \right.
\end{displaymath}
Note that in case $G=\Z/2\Z$ and $M=\f$, the cohomology of $\Z/2\Z$ with coefficients in $\f$ equals
$$H^n(G,\f)=\f  \textrm{~~ for ~~ }n\geq 0.$$
\end{ex}

\subsection{Definition}
% Pour la definition, on a besoin aussi de X non compacte pour la suite d'une paire.
Let $G$ be a group, $A$ a commutative ring, $M$ an $A$-module with a $G$-action and let $X$ be a topological space with a $G$-action. Then, for $n\in \Z$, we define the $n$-th equivariant Borel-Moore homology group with coefficients in $M$ by:
$$H_n(X;G,M)=R^{-n}\Hom _G(R\Gamma_c A,M).$$

Note that in the case of a point, one recovers the cohomology of the group with coefficients in $M$: $H_n(\textbf{pt};G,M)=H^{-n}(G,M).$

\begin{rmk}\label{cell} One can give an alternative definition using cellular chain complexes (see \cite{Derval}).
Let $C_{\cdot}(X)$ be a CW-complex with a $G$-action on the cells. Let $F_{\cdot}$ be a projective resolution of $\Z$ via $\Z[G]$-modules. Then the $n$-th equivariant Borel-Moore homology group of $X$ with coefficients in $M$ is equal to
$$H_n(X;G,M)=H^{-n}(\Hom _G(F_{\cdot},C_{\cdot} \otimes M)).$$
The main interest of such a description is to enable computations by using cell decompositions.
\end{rmk}

\begin{ex}\label{sphere} Consider an action of $G=\Z/2\Z$ on the sphere $S^d$ with at least one fixed point. Take $A=\f$. Then one can choose, as a cell decomposition, the one composed of a fixed point $\textbf{pt}$ as a 0-cell, and $S^d \setminus \textbf{pt}$ as a $d$-cell. Then the action of $G$ on the cell decomposition does not depend on the action of $G$ on the sphere (apart from the possible existence of a fixed point).
\end{ex}

\begin{rmk}\begin{flushleft}\end{flushleft}
\begin{enumerate}
\item In the same way, one can define equivariant cohomology groups with coefficients in $M$ by
$$H^n(X;G,M)=R^n \Hom _G(R\Gamma A,M).$$
It admits also an analogous description in terms of cell decompositions.
\item In case $G$ is trivial, the equivariant Borel-Moore homology groups and the equivariant cohomology groups coincide with the non-equivariant ones.
\item Note that this definition of equivariant homology, as a mixture of cohomology of groups and homology of varieties, is not the classical one. However, this definition is the natural dual to equivariant cohomology (see \cite{VH}). In particular, Poincar\'e duality holds for $M=\f$.
\end{enumerate}
\end{rmk}
%%%%%%%%%%%%%%%%%%%%%%%%%%%%%%%%%%%%%%%%%%%%%%%%%%%%%%%%%%%%%%%%%%%%%%%%%%
\subsection{Some properties}
% suite exacte d'une paire, dualite poincare, action triviale
In this section, we recall some properties of equivariant homology without any proof. We refer to \cite{VH}, \cite{Derval} for more details.
\subsubsection{Long exact sequence of a pair}
 Let $i:Z \hookrightarrow X$ denote the inclusion of a closed $G$-invariant subspace $Z$ into a $G$-space $X$. Denote by  $j:U=X \setminus Z \hookrightarrow X$ the inclusion of the complement. Then the sequence
$$\xymatrix{ \cdots \ar[r] & H_n(Z;G,M) \ar[r]^{i_*} & H_n(X;G,M) \ar[r]^{j^*}  & H_n(U;G,M) \ar[r] & H_{n-1}(Z;G,M) \ar[r]  & \cdots }$$
is exact.

\subsubsection{Trivial action}\label{triv}
In case $A$ is a field and the actions of $G$ on $X$ and on $M$ are trivial, one recovers the equivariant homology groups of a locally compact space $X$ of finite cohomological dimension from the homology of $X$ and the cohomology of $G$. More precisely, the cap product induces an isomorphism:
$$\oplus _{p-q=n} \big(H_p(X,M) \otimes H^q(G,A) \big) \dans H_n(X;G,M).$$
\vskip 3mm

\begin{ex}\label{sphere-triv} Let us take $G=\Z/2\Z$ and $M=A=\f$. The $G$-equivariant homology of $S^d$ with trivial action is
\begin{displaymath}
H_n(S^d,G,\f) = \left\{ \begin{array}{lll}
\f   & \textrm{if } 1 \leq n \leq d,\\
(\f)^2   & \textrm{if } n \leq 0\\
 0    & \textrm{if } n > d.
\end{array} \right.
\end{displaymath}
In particular, it follows from example \ref{sphere} that these groups are the equivariant homology groups of $S^d$ for any $G$-action which is not free.
\end{ex}
%%%%%%%%%%%%%%%%%%%%%%%%%%%%%%%%%%%%%%%%%%%%%%%%%%%%%%%%%%%%%%%%%%%%%%%%%%
\subsubsection{Fixed point free action}\label{free}
If the action is free on $X$ and trivial on $M$, the equivariant cohomology is directly related to the homology of the quotient via the isomorphism:
$$H_n(X;G,M)\simeq H_n(X/G;M)$$
for any $n\in \Z$.
%%%%%%%%%%%%%%%%%%%%%%%%%%%%%%%%%%%%%%%%%%%%%%%%%%%%%%%%%%%%%%%%%%%%%%%%%%
\subsubsection{K\"unneth formula}
% enoncer la formule de Kunneth, on en a besoin pour une demo
Let $G,G'$ be finite groups and $X,Y$ be topological spaces of finite cohomological dimension with respectively a $G$-action and a $G'$-action. Let $k$ be a field with trivial $G$-action and $G'$-action. Then, for $n\in \Z$, we have isomorphisms:
$$H_n(X\times Y;G\times G',k) \simeq \oplus _{p+q=n} \big( H_p(X;G,k) \otimes H_q(Y;G',k) \big ).$$ 

%%%%%%%%%%%%%%%%%%%%%%%%%%%%%%%%%%%%%%%%%%%%%%%%%%%%%%%%%%%%%%%%%%%%%%%%%%
\subsubsection{Poincar\'e duality}
In order to avoid orientation issues, let us fix $M=A=\f$.
One can define cup and cap product in the equivariant setting. Moreover, for $X$ a compact $d$-dimensional manifold, one can define the $G$-equivariant fundamental class $\mu_X \in H_n(X;G,\f)$ of $X$. Then cap product with $\mu_X$ induces an isomorphism:
$$\xymatrix{H^n(X;G,\f) \ar[r]^{\cap \mu_X} & H_{d-n}(X;G,\f).}$$

%%%%%%%%%%%%%%%%%%%%%%%%%%%%%%%%%%%%%%%%%%%%%%%%%%%%%%%%%%%%%%%%%%%%%%%%%%
\subsubsection{Hochschild-Serre spectral sequence}
% et exemples de calculs: cercles, sphere avec un point fixe
Equivariant (co)homology comes with spectral sequences that make computations possible (\cite{groth}). We present below the Hochschild-Serre spectral sequence.
\begin{prop} Let $G$ be a group, $X$ a locally compact space of finite cohomological dimension with a $G$-action and $M$ an $A$-module with a $G$-action. Then there exists a spectral sequence
$$E^2_{p,q}=H^{-p}(G,H_q(X,M)) \Rightarrow H_{p+q}(X;G,M)$$
functorial in $X$ and $M$.
\end{prop}

In the particular case $G=\Z/2\Z$, the Hochschild-Serre spectral sequence is greatly simplified:
\begin{displaymath}
E^2_{p,q} = \left\{ \begin{array}{lll}
 H_q(X,M)^G  & \textrm{if~~}p=0,\\
 H^{1}(G,H_q(X,M))  & \textrm{if~~}p\leq -1.\\
\end{array} \right.
\end{displaymath}

Moreover, if $M=A=\f$, the following results guarantee the vanishing of some differentials under the condition of existence of a fixed point.
\begin{lemma}\label{ss} Let $X$ be a smooth compact connected space of finite dimension $d$ endowed with an action of $G=\Z/2\Z$. The differentials of the Hochschild-Serre spectral sequence
$$E^2_{p,q}=H^{-p}(G,H_q(X,\f)) \Rightarrow H_{p+q}(X;G,\f)$$
\begin{enumerate}
\item having source $E^r_{p,0}$ are trivial for any $p \leq 0,~r \geq2$ if and only if $X^G \neq \emptyset$,
\item having target $E^r_{p,d}$ are trivial for any $p \leq 0,~r \geq2$ if and only if $X^G \neq \emptyset$.
\end{enumerate}
\end{lemma}

\begin{ex}\label{ss-sphere}\begin{flushleft}\end{flushleft}
\begin{enumerate}
\item \label{fixe} Consider a $d$-sphere, $d\geq 1$, with an action of $G=\Z/2\Z$ that admits at least a fixed point. The $G$-action on the homology of $S^d$ is trivial since the homology groups of $S^d$ are either $\f$ or $\{0\}$, therefore $E^2_{p,q}=\f$ if $p\leq 0$ and $q=0,d$ whereas $E^2_{p,q}=0$ otherwise. 

In the diagram, only two half-lines are non-zero.

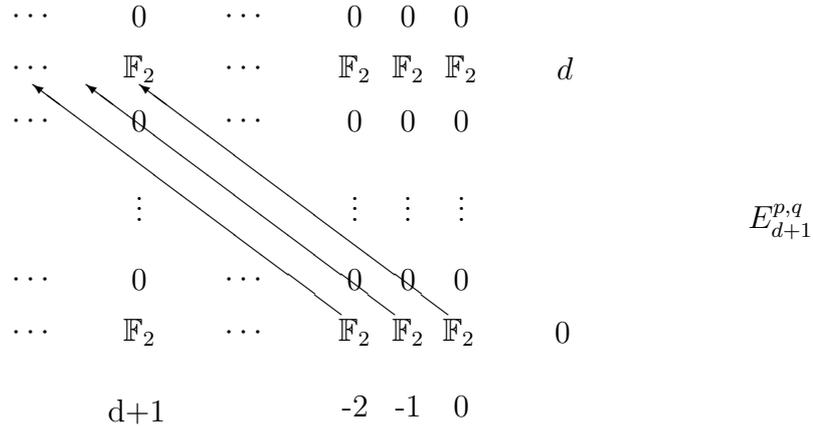
\unitlength=.7mm
\begin{picture}(130,100)(-40,-10)
%\put(-20,0){\vector(1,0){85}}
\put(58,3){\vector(-4,3){58}}
\put(48,3){\vector(-4,3){58}}
\put(38,3){\vector(-4,3){58}}
\put(0,-15){\makebox(-1,-1){d+1}}
%\put(70,70){\makebox(5,5){q}}
\put(80,0){\makebox(-1,-1){0}}

\put(0,0){\makebox(0,-1){$\f$}}
\put(60,0){\makebox(0,-1){$\f$}}
\put(20,0){\makebox(0,-1){$\cdots$}}\put(-20,0){\makebox(0,-1){$\cdots$}}
\put(50,0){\makebox(1,-1){$\f$}}\put(50,-15){\makebox(1,1){-1}}
\put(40,0){\makebox(1,-1){$\f$}}\put(40,-15){\makebox(1,1){-2}}

\put(0,60){\makebox(0,-1){$0$}}
\put(20,60){\makebox(0,-1){$\cdots$}}\put(-20,60){\makebox(0,-1){$\cdots$}}
\put(40,60){\makebox(1,-1){$0$}}
\put(60,60){\makebox(1,-1){$0$}}
\put(50,60){\makebox(1,-1){$0$}}

\put(0,10){\makebox(0,-1){$0$}}
\put(20,10){\makebox(0,-1){$\cdots$}}\put(-20,10){\makebox(0,-1){$\cdots$}}
\put(40,10){\makebox(1,-1){$0$}}
\put(60,10){\makebox(1,-1){$0$}}
\put(50,10){\makebox(1,-1){$0$}}

\put(0,25){\makebox(0,-1){$\vdots$}}
%\put(20,25){\makebox(0,-1){$\cdots$}}
\put(40,25){\makebox(1,-1){$\vdots$}}
\put(60,25){\makebox(1,-1){$\vdots$}}
\put(50,25){\makebox(1,-1){$\vdots$}}

\put(0,50){\makebox(0,-1){$\f$}}
\put(20,50){\makebox(0,-1){$\cdots$}}\put(-20,50){\makebox(0,-1){$\cdots$}}
\put(40,50){\makebox(1,-1){$\f$}}
\put(60,50){\makebox(1,-1){$\f$}}
\put(50,50){\makebox(1,-1){$\f$}}

\put(0,40){\makebox(0,-1){$0$}}
\put(20,40){\makebox(0,-1){$\cdots$}}\put(-20,40){\makebox(0,-1){$\cdots$}}
\put(40,40){\makebox(1,-1){$0$}}
\put(60,40){\makebox(1,-1){$0$}}
\put(50,40){\makebox(1,-1){$0$}}

\put(80,50){\makebox(0,-1){$d$}}
\put(60,-15){\makebox(1,1){$0$}}

\put(120,20){\makebox(1,1){$E_{d+1}^{p,q}$}}
\end{picture}
\vskip 10mm
Now, the differentials, up to the level 2, all vanish thanks to lemma \ref{ss}, therefore the Hochschild-Serre spectral sequence degenerates at the level $E^2$. As a consequence:
\begin{displaymath}
H_n(S^d,G,\f) = \left\{ \begin{array}{lll}
\f   & \textrm{if } 1 \leq n \leq d,\\
(\f)^2   & \textrm{if } n \leq 0\\
 0    & \textrm{if } n \geq d+1.
\end{array} \right.
\end{displaymath}
Note that we recover here the previous results obtained in examples \ref{sphere} and \ref{sphere-triv}.
\item\label{central} Consider the unit sphere $S^1$ with the action of $G=\Z/2\Z$ given by the central symmetry. There is no fixed point, thus at least one of the differentials having source $E^2_{p,0}$ is not trivial by lemma \ref{ss}. Actually all of these differentials are the same, therefore the Hochschild-Serre spectral sequence degenerates at the level $E^3$ and
\begin{displaymath}
H_n(S^1,G,\f) = \left\{ \begin{array}{ll}
\f   & \textrm{if } n \in \{0,1\},\\
0   & \textrm{otherwise.}
\end{array} \right.
\end{displaymath}
In particular we recover that $H_n(S^1,G,\f)$ is isomorphic to $H_n(S^1/G,\f)$ since the action is free (cf. section \ref{free}).
\end{enumerate}
\end{ex}

%%%%%%%%%%%%%%%%%%%%%%%%%%%%%%%%%%%%%%%%%%%%%%%%%%%%%%%%%%%%%%%%%%%%%%%%%%
\subsubsection{ Negative degree in the case $G=\Z/2\Z$}\label{negative}
Take $n<0$. The equivariant homology group $H_n(X;G,\f)$ can be computed in terms of the homology of the fixed point set $X^G$ as follows:
$$H_n(X;G,\f) \simeq \oplus_{i \geq 0}H_i(X^G;\f).$$

%%%%%%%%%%%%%%%%%%%%%%%%%%%%%%%%%%%%%%%%%%%%%%%%%%%%%%%%%%%%%%%%%%%%%%%%%%
%%%%%%%%%%%%%%%%%%%%%%%%%%%%%%%%%%%%%%%%%%%%%%%%%%%%%%%%%%%%%%%%%%%%%%%%%%
%%%%%%%%%%%%%%%%%%%%%%%%%%%%%%%%%%%%%%%%%%%%%%%%%%%%%%%%%%%%%%%%%%%%%%%%%%
\section{Equivariant virtual Betti numbers}\label{sect-eq-virt}

In this part, we introduce equivariant virtual Betti numbers for arc-symmetric sets, with action of a finite group $G$. We put emphasis on the case $G=\Z/2\Z$ in the examples.

%%%%%%%%%%%%%%%%%%%%%%%%%%%%%%%%%%%%%%%%%%%%%%%%%%%%%%%%%%%%%%%%%%%%%%%%%%
\subsection{Arc-symmetric sets}
Arc-symmetric sets have been introduced by K. Kurdyka \cite{KK1} in
order to study ``rigid components'' of real algebraic varieties. The category of arc-symmetric sets contains real algebraic
varieties and, in some sense, this category has better behaviour
than the category of real algebraic varieties, maybe closer to complex
algebraic varieties (see remark \ref{quotient-alg} below for instance). For a detailed treatment of arc-symmetric sets, we refer
to \cite{KK1,fichou}. Nevertheless, let us make precise the definition of these sets.

We fix a compactification of $\mathbb R^n$, for instance $\mathbb R^n
\subset \mathbb P^n(\mathbb R)$.
\begin{defi} Let $X \subset \mathbb P^n$ be a semi-algebraic set. We say
  that $X$ is arc-symmetric if,
 for every real analytic arc $\gamma :
  ]-1,1[ \longrightarrow  \mathbb P^n$ such that $\gamma
  (]-1,0[) \subset X$, there exists $\epsilon  > 0$ such that $\gamma
  (]0,\epsilon [) \subset X$.
\end{defi}

One can think of arc-symmetric sets as the
smallest category, denoted $\mathcal {AS}$, stable under boolean operations and containing compact
real algebraic varieties and their connected components.
Note that by the dimension of an arc-symmetric set, we mean its dimension as a semi-algebraic set. 

Arc-symmetric sets are related to connected components of real algebraic varieties. The definition of irreducibility for arc-symmetric sets is the classical one.

\begin{prop}(\cite{KK1})\label{ecl} Let $X \in \mathcal {AS}$ be irreducible. Let
  $Z$ be a compact
  real algebraic variety containing $X$ with $\dim Z=\dim X$, and $\pi :
  \widetilde Z \longrightarrow Z$ a resolution of singularities for
  $Z$. Denote by $\Reg(X)$ the complement in $X$ of the singular points of $Z$.
Then there exists a unique connected component $ \widetilde X$ of $
  \widetilde Z$ such that $ \pi (\widetilde X)$ is equal to the euclidean closure $\overline {\Reg(X)}$ of $\Reg(X)$.
\end{prop}

In particular, the following corollary specifies what nonsingular
arc-symmetric sets look like.
Note that by an isomorphism between arc-symmetric sets, we mean a birational
map containing the arc-symmetric sets in its support. Moreover, a
nonsingular irreducible arc-symmetric set is an irreducible arc-symmetric set whose intersection
with the singular locus of its Zariski closure is empty.

\begin{cor} Compact nonsingular arc-symmetric sets are isomorphic to
  unions of
  connected components of compact nonsingular real algebraic varieties.
\end{cor}

%%%%%%%%%%%%%%%%%%%%%%%%%%%%%%%%%%%%%%%%%%%%%%%%%%%%%%%%%%%%%%%%%%%%%%%%%%
\subsubsection{Equivariant Grothendieck group of arc-symmetric sets}\label{}
Let $G$ be a finite group. The equivariant Grothendieck group $K^0(\AS,G)$ of arc-symmetric sets is defined to be the group generated by the symbols $[X,G]$ for $X$ a $G$-equivariant isomorphism class of arc-symmetric sets endowed with an algebraic $G$-action, which means an action of $G$ on $X$ that comes from an algebraic action on the Zariski closure of $X$, subject to the relations
\begin{itemize}
\item $[X,G]=[Y,G]-[X \setminus Y,G]$ if $Y$ is a closed
  $G$-invariant arc-symmetric subset of $X$,
\item $[V,G]$=$[\mathbb A^n \times X,G]$ if $V \longrightarrow X$ is the restriction to the arc-symmetric set $X$ of a vector bundle on the Zariski closure ${\overline X}^Z$ of $X$, with a linear $G$-action over the action on ${\overline X}^Z$. Here, on the right hand side, the action of $G$ on $\mathbb A^n$ is trivial.
\end{itemize}
Note that by definition of the $G$-action, an arc-symmetric set admits an equivariant compactification $X\hookrightarrow \overline X$, which means there exists a compact arc-symmetric set $\overline X$ containing $X$, endowed with a $G$-action whose restriction to $X$ coincides with the $G$-action on $X$.
Note also that the product of varieties, with the diagonal action, induces a commutative ring structure on $K^0(\AS,G)$.

There exists a simple presentation of this rather complicated group in terms of connected components of compact nonsingular real algebraic varieties and blowings-up (cf. \cite{FB, fichou}).

For a connected component $X$ of a compact nonsingular real algebraic variety $Z$, and $C \subset Z$ a closed nonsingular subvariety, we will denote by $\Bl_{C \cap X} X$ the restriction to $X$ of the blowing-up $\Bl_C Z$ of $Z$ along $C$, called the blowing-up of $X$ along $C \cap X$.

\begin{thm}\label{bitt}
The group $K^0(\AS,G)$ is the abelian group generated by the isomorphism classes of connected components of compact nonsingular affine real algebraic varieties with $G$-action, subject to the relations $[\emptyset,G]=0$ and $[\Bl_{C\cap X} X,G]-[E,G]=[X,G]-[C \cap X,G]$, where $X$ is a connected component of a compact nonsingular real algebraic variety $Z$ with $G$-action, $C\subset Z$ is a closed nonsingular $G$-invariant subvariety, $\Bl_{C\cap X} X$ is the blowing-up of $X$ along $C \cap X$ and $E$ is the exceptional divisor of this blowing-up.
\end{thm}

\begin{rmk}\label{proof-bitt} One can prove the theorem by induction on the dimension of the varieties. The principle of the proof is to express the class of
\begin{enumerate}
\item a nonsingular arc-symmetric set $X$ in terms of that of an equivariant nonsingular compactification $X \hookrightarrow \overline X$ (that exists by \cite{hiro}) and of its complement $\overline X \setminus X$ that has dimension strictly less than that of $X$. It can be proven that the class $[X]=[\overline X]-[\overline X \setminus X]$ does not depend on the choice of the compactification via the equivariant weak factorisation theorem \cite{weak},
\item a singular arc-symmetric set in terms of the class of nonsingular arc-symmetric sets via a stratification with nonsingular strata.
\end{enumerate}
\end{rmk}

%%%%%%%%%%%%%%%%%%%%%%%%%%%%%%%%%%%%%%%%%%%%%%%%%%%%%%%%%%%%%%%%%%%%%%%%%%
\subsection{Equivariant virtual Betti numbers}\label{sect-virt}
Let us denote by $b_i^G(X)$ the i-$th$ $G$-equivariant Betti number of $X$, that is the dimension over $\f$ of the i-$th$ equivariant homology group with coefficient in $\f$ of a compact nonsingular arc-symmetric set with a $G$-action.

We begin with a useful lemma.

\begin{lemma}\label{add} Let $X$ be a connected component of a compact nonsingular real algebraic variety $Z$ with $G$-action and $C\subset Z$ be a closed nonsingular $G$-invariant subvariety. Denote by $\Bl_{C\cap X} X$ the blowing-up of $X$ along $C \cap X$ and by $E$ the exceptional divisor of this blowing-up. Then the sequence 
$$ 0 \dans H_i(E;G,\f)   \dans  H_i(\Bl_{C\cap X} X;G,\f)\oplus H_i(C\cap X;G,\f)  \dans H_i(X;G,\f)  \dans 0$$
is exact for $i\in \Z$.
\end{lemma}

As in the non-equivariant case (\cite{MCP}), this is a consequence of Poincar\'e duality and the existence of the long exact sequence of a pair for equivariant homology. We recall briefly the proof for the convenience of the reader.

\begin{proof} For sake of simplicity, we replace $C \cap X$ by $C$ in the proof.
Note that $\pi_i: H_i(\Bl_C X;G,\f) \dans H_i(X;G,\f) $ is surjective by Poincar\'e duality. Indeed, the diagram 
$$\xymatrix{H^{d-i}(X;G,\f) \ar[rr]^{\cap \mu_X}_{\simeq}\ar[d]_{\pi^{d-i}}& & H_i(X;G,\f)\\
           H^{d-i}(\Bl_C X;G,\f) \ar[rr]_{\cap \mu_{\Bl_C X}}^{\simeq}            & & H_i(\Bl_C X;G,\f) \ar[u]_{\pi _i}}$$
is commutative, where $d=\dim X$ and $\pi:\Bl_C X \dans X$ is the blowing-up.

Moreover, the long exact sequences of the pairs $(X,C)$ and $(\Bl_C X,E)$ give rise to a commutative diagram
$$\xymatrix{\cdots \ar[r] & H_{i-1}(X\setminus C;G,\f)  \ar[r]  & H_{i}(E;G,\f)  \ar[r]   & H_i (X;G,\f) \ar[r] &\cdots \\
            \cdots \ar[r] & H_{i-1}((\Bl_C X) \setminus E;G,\f) \ar[r]\ar[u]^{\simeq} & H_i(C;G,\f) \ar[r]\ar[u]  & H_i(\Bl_C X;G,\f) \ar[r]\ar[u] &\cdots} $$
where the vertical arrows are induced by $\pi$. 
We obtain now the announced exact sequence by straightforward diagram chasing.
\end{proof}

The following theorem is the main result of this paper.

\begin{thm}\label{main} Let $G$ be a finite group.
For each integer $i \in \mathbb Z$ there exists a unique group homomorphism $\beta_i^G(.):K^0(\AS,G) \dans \Z$ such that $\beta_i^G(X)=b_i^G(X)$ for $X$ compact and nonsingular.
\end{thm}

\begin{proof} The construction, invariance and additivity are direct consequences of theorem \ref{bitt} and lemma \ref{add} and remain valid even for non-finite groups. The fact that $\beta_i^G(V)=\beta_i^G(\mathbb A^n \times X)$ for a vector bundle $V$ with a linear $G$-action over the action of a finite group $G$ on $X$ will follow from proposition \ref{mult} below.
\end{proof}

We call $\beta_i^G (X)$ the $i$-th equivariant virtual Betti number of $X$.
One may define the equivariant virtual Poincar\'e series of an arc-symmetric set with a $G$-action to be the formal power series
$$\beta^G(X)=\sum_{i \in \Z} \beta_i^G(X) u^i \in \Z[[u,u^{-1}]].$$
Denote also by $b^G(X)$ the equivariant Poincar\'e series of a compact nonsingular arc-symmetric set $X$ with a $G$-action.

\begin{rmk} For $G$ the trivial group, the equivariant virtual Poincar\'e series $\beta^G$ is equal to the virtual Poincar\'e polynomial $\beta$ (cf. \cite{MCP, fichou}).
\end{rmk}

Similarly to the case of the virtual Poincar\'e polynomial (cf. \cite{fichou}), the equivariant virtual Poincar\'e series remains invariant not only under algebraic isomorphisms between arc-symmetric sets, but also under Nash isomorphisms. 
A Nash isomorphism
between arc-symmetric sets $X_1,X_2 \in \mathcal{AS}$ is the
restriction of an analytic
and semi-algebraic isomorphism between compact semi-algebraic and real analytic sets $Y_1,Y_2$
containing $X_1,X_2$ respectively. 

The proof in the non-equivariant case carries over to the equivariant case.

\begin{thm}\label{nash-iso} Let $X_1,X_2 \in \mathcal{AS}$ be endowed with a $G$-action. If there exists a Nash isomorphism between $X_1$ and $X_2$, equivariant with respect to the $G$-actions, then the equivariant virtual Poincar\'e series of $X_1$ and $X_2$ coincide.
\end{thm}

An important result about the equivariant virtual Poincar\'e series of an arc-symmetric set is that it keeps its dimension. In particular the equivariant virtual Poincar\'e series takes its values in $\Z[[u^{-1}]][u]$.

\begin{prop} The formal power series $\beta^G(X)$ of an arc-symmetric set $X$ with a $G$-action has finite degree equal to the dimension of $X$.
\end{prop}
\begin{proof} We prove the result by induction on the dimension of $X$. Let us divide the proof into several steps.
\begin{enumerate}
\item\label{p1} In case $X$ is compact and nonsingular, the equivariant virtual Betti numbers coincide with the equivariant Betti numbers. Therefore the result is true because the $n$-th equivariant homology group $H_n(X;G,\f)$ is equal to zero for $n$ strictly greater than the dimension of $X$ whereas it is nonzero when $n$ equals the dimension of $X$.

\item To initiate the induction, note that in dimension zero, the arc-symmetric sets are compact and nonsingular. So the result is true by (\ref{p1}). 
\item\label{p3} Now, assume $X$ is a nonsingular arc-symmetric set of dimension $d$, possibly not compact. Choose a nonsingular equivariant compactification $X \hookrightarrow \overline X$ with $\dim \overline X \setminus X < \dim X$ (this is possible by \cite{hiro}). Then
$$\beta^G(X)=\beta^G(\overline X)-\beta^G(\overline X \setminus X)$$
by additivity of the equivariant virtual Poincar\'e series. So the degree of $\beta^G(X)$ is finite by (\ref{p1}) and by the inductive assumption.

\item Finally, in the case of a singular $X$, we achieve the proof by stratifying $X$ into nonsingular strata for which the result is true by (\ref{p3}). 
\end{enumerate}
\end{proof}

\begin{rmk}\label{point} However the equivariant virtual Poincar\'e series does not admit a finite valuation in general. It particular for $G=\Z/2\Z$ and $X$ a single point:
$$\beta^G(\textbf{pt})= \sum_{i \leq 0} u^i=\frac{u}{u-1}.$$
\end{rmk}

\begin{ex}\label{ex-beta} In the following examples $G=\Z/2\Z$.\begin{flushleft}\end{flushleft}
\begin{enumerate}
\item\label{point-ech} The equivariant virtual Poincar\'e series of two points $\{p,q\}$ that are inverted by the action is
$$\beta^G(\{p,q\})=b^G(\{p,q\})=1$$
because such a variety is compact and nonsingular.
\item In the same way, if we consider the $d$-sphere with the action given by the central symmetry, then:
$$\beta^G(S^d)=b^G(S^d)=1+u^d$$
thanks to example \ref{ss-sphere}.\ref{fixe}, whereas when we consider it with a $G$-action that admits a fixed point:
$$\beta^G(S^d)=b^G(S^d)=u^d+\cdots+u+2\frac{u}{u-1}$$
as we computed in example \ref{ss-sphere}.\ref{central}.
\item By additivity, the equivariant virtual Poincar\'e series of the real affine line endowed with the trivial $G$-action is equal to:
$$\beta^G(\mathbb A^1)=\beta^G(S^1)-\beta^G(\textbf{pt})=\sum_{i \leq 1} u^i=\frac{u^2}{u-1},$$
by example \ref{sphere-triv}.
\item\label{ex-beta-4} Now, consider the real affine line $\mathbb A^1$ endowed with any $G$-action. Again 
$$\beta^G(\mathbb A^1)=\frac{u^2}{u-1}.$$
Actually the $G$-action on the compactification $A^1 \hookrightarrow S^1$ admits at least a fixed point, so the result follows from example \ref{ss-sphere}.\ref{fixe}. Similarly:
$$\beta^G(\mathbb A^d)=\frac{u^{d+1}}{u-1}.$$
\end{enumerate}
\end{ex}

Let us notice that the result of example \ref{ex-beta}.\ref{ex-beta-4} can be rewritten in the form
$$\beta^G(\mathbb A^d)=u^d \beta^G(\textbf{pt}).$$
The following fact consists in a generalisation of this remark. Moreover, it enables us to complete the proof of theorem \ref{main}.

% pour la formule de DL
\begin{prop}\label{mult} Let $X$ be an arc-symmetric set with a $G$-action, and consider the affine variety $\mathbb A^d$ with any $G$-action. Then
$$\beta^G(X \times \mathbb A^d)=u^d \beta^G(X)$$
where the left hand side is endowed with the diagonal $G$-action.
\end{prop}

\begin{proof} To begin with, let us tackle the compact nonsingular case. Take $X$ compact and nonsingular, and compactify $X \times \mathbb A^d$ into $X \times \mathbb S^d$. Choose as a $G$-invariant cell decomposition of $X \times \mathbb S^d$ the product of any $G$-invariant cell decomposition of $X$ and of the cell decomposition of $S^d$ made of $\mathbb A^d$ as a $d$-cell and of $\mathbb S^d \setminus\mathbb A^d$ as a $0$-cell. Now, the action of $G$ on this cell decomposition of $X \times \mathbb S^d$ is the same as the action of $G \times \{1\}$ on $X \times \mathbb S^d$ since the $G$-action on the cell decomposition of $\mathbb S^d$ is trivial. Therefore 
$$b^G(X \times S^d)=b^{G\times \{1\}}(X \times S^d).$$
Now, by the K\"unneth formula:
$$b^{G\times \{1\}}(X \times S^d)=b^G(X) b^{\{1\}}(S^d)=(1+u^d) b^G(X),$$
and so 
$$b^G(X \times \mathbb A^d)=b^G(X \times S^d)-b^G(X \times \textbf{pt})=(1+u^d) b^G(X)- b^G(X)=u^d  b^G(X).$$

Now, let us deal with the general case. We prove it by induction on the dimension of $X$. We divide the proof into several steps.
\begin{enumerate}
\item In dimension zero, the result is true since $0$-dimensional arc-symmetric sets are compact and nonsingular.
\item\label{cas-lisse} Let $X$ be a $n$-dimensional nonsingular arc-symmetric set endowed with a $G$-action. Compactify $X$ into a nonsingular variety $X \hookrightarrow \overline X$ with a $G$-action compatible with that of $X$, and such that the dimension of the complement $\overline X \setminus X$ is strictly less that $n$. Then compactify $X\times \mathbb A^d$ into $X\times \mathbb A^d \hookrightarrow \overline X \times S^d$. Note that the $G$-action on $S^d$ admits a fixed point. So
$$\overline X \times S^d= X\times \mathbb A^d \cup X \times \textbf{pt}  \cup  (\overline X \setminus X )\times S^d$$
and thus
$$\beta^G(X\times \mathbb A^d )=\beta^G(\overline X \times S^d)-\beta^G(X \times \textbf{pt})-\beta^G((\overline X \setminus X )\times S^d).$$
However 
$$\beta^G(\overline X \times S^d)=(1+u^d) \beta^G(\overline X)$$
 because $\overline X $ is compact and nonsingular, and moreover
$$\beta^G((\overline X \setminus X )\times S^d)=(1+u^d) \beta^G(\overline X \setminus X )$$
by the inductive assumption. Therefore, by additivity of the equivariant virtual Poincar\'e series:
$$\beta^G(X\times \mathbb A^d )=(1+u^d)\big(\beta^G(\overline X)-\beta^G(\overline X \setminus X )\big)-\beta^G(X)$$
$$=(1+u^d)\beta^G(X)-\beta^G(X)=u^d \beta^G(X).$$
\item Finally, for $X$ general, stratify $X$ into nonsingular strata and apply (\ref{cas-lisse}).
\end{enumerate}
\end{proof}
%%%%%%%%%%%%%%%%%%%%%%%%%%%%%%%%%%%%%%%%%%%%%%%%%%%%%%%%%%%%%%%%%%%%%%%%%%
\subsection{Particular actions}
%%%%%%%%%%%%%%%%%%%%%%%%%%%%%%%%%%%%%%%%%%%%%%%%%%%%%%%%%%%%%%%%%%%%%%%%%%
\subsubsection{Trivial actions}
In case the action is trivial, one can relate the equivariant virtual Poincar\'e series to the virtual Poincar\'e polynomial.
\begin{prop}\label{action-triv} If $X$ is an arc-symmetric set endowed with a trivial $G$-action, then
$$\beta^G(X)=\beta(X)\beta^G(\textbf{pt}).$$
\end{prop}

\begin{proof} Remark first that the result is true if $X$ is compact and nonsingular (cf. section \ref{triv}). Now, let us prove the proposition by induction on the dimension of $X$. In dimension zero, the result is true by additivity. Now, take $X$ endowed with a trivial $G$-action. One can assume $X$ to be nonsingular, otherwise stratify $X$ into nonsingular strata and use the additivity to reduce to this case. Then, choose a nonsingular compactification $\overline X$ of $X$ such that the dimension of $\overline X \setminus X$ is strictly less than that of $X$, and endow it with the trivial action of $G$. Therefore we obtain an equivariant compactification of $X$ and: 
$$\beta^G(X)=\beta^G(\overline X)-\beta^G(\overline X \setminus X)$$
by additivity. Now by the induction assumption we get: $$\beta^G(\overline X \setminus X)=\beta(\overline X \setminus X)\beta^G(\textbf{pt}).$$  Thanks to the compact nonsingular case, we obtain: $$\beta^G(\overline X)=\beta(\overline X)\beta^G(\textbf{pt}).$$ We conclude by additivity of the virtual Poincar\'e polynomial.
\end{proof}
%%%%%%%%%%%%%%%%%%%%%%%%%%%%%%%%%%%%%%%%%%%%%%%%%%%%%%%%%%%%%%%%%%%%%%%%%%
\subsubsection{Free actions}
The equivariant homology of a topological space endowed with a free action is equal to the cohomology of its quotient (section \ref{free}). The next result states the analog for the equivariant virtual Poincar\'e series in the compact case. In particular, we prove that the quotient, viewed as a semi-algebraic subset of the quotient of the Zariski closure (which is itself a semi-algebraic set in general), is still an arc-symmetric set in this case.
\begin{prop}\label{action-free} Let $X$ be a compact arc-symmetric set endowed with a free $G$-action. Then the quotient $X/G$ is an arc-symmetric set and
$$\beta^G(X)=\beta(X/G).$$
\end{prop}

\begin{proof} Let us prove first that the quotient is arc-symmetric. First, the quotient $X/G$ remains semi-algebraic because so is $X$. Assume that $X/G$ is a semi-algebraic set in $\mathbb R^n\subset \mathbb P^n(\mathbb R)$. Take $\gamma :
  ]-1,1[ \longrightarrow  \mathbb P^n$ an analytic arc such that $\gamma
  (]-1,0[) \subset X/G$. Then $\gamma (0)\in X/G$ because $X$ is closed and $\gamma$ admits an analytic lifting $\widetilde \gamma :]-1,0] \longrightarrow X$ since the action is free on $X$. But $X$ is arc-symmetric, therefore there exists $\epsilon  > 0$ such that $\widetilde \gamma$ can be extended analytically to $]-1,\epsilon [$ with
  $\widetilde \gamma (]0,\epsilon [) \subset X$. However the push-forward of $\widetilde \gamma$ by $\pi$ must coincide with $\gamma$ on $]0,\epsilon [$, thus $\gamma (]0,\epsilon [) \subset X/G$.

For the second part, we prove the result, once more, by induction on the dimension. Let $\pi:Y \dans X$ be an equivariant resolution of the singularities of $X$. Then $Y$ is also endowed with a $G$-action which must be free because so is the $G$-action on $X$. Denote by $E$ the exceptional divisor and by $C$ the centre of the resolution. Then $E$ and $C$ have a $G$-action and $Y \setminus E$ is isomorphic to $X \setminus C$ via $\pi$.

Now $\beta^G(E)=\beta(E/G)$ and $\beta^G(C)=\beta(C/G)$ by the induction assumption, and $\beta^G(Y)=\beta(Y/G)$ since $Y$ is compact and nonsingular. Therefore
$$\beta^G(X)=\beta(Y/G)+\beta(C/G)-\beta(E/G).$$
On the other hand, the quotient $(Y \setminus E)/G$ is isomorphic to $(X \setminus C)/G$ so $$\beta((Y \setminus E)/G)=\beta((X \setminus C)/G).$$
Then by additivity 
$$\beta(X/G)=\beta((X \setminus C)/G)+\beta(C/G)$$
and
$$\beta(Y/G)=\beta((Y \setminus E)/G)+\beta(E/G).$$
Finally, once more by additivity
$$\beta(X/G)=\beta((X \setminus C)/G)+\beta(C/G)=\beta((Y \setminus E)/G)+\beta(C/G)=\beta(Y/G)-\beta(E/G)+\beta(C/G).$$
\end{proof}

\begin{rmk}\label{quotient-alg} If $X$ is no longer compact or if the action is no longer free, then the quotient space is no longer an arc-symmetric set in general. For instance, a reflexion on a circle gives as quotient a half-circle. Moreover, even if $X$ is a compact real algebraic variety, then the quotient space is not algebraic in general but only arc-symmetric. For instance, consider the real algebraic plane curve $C$ given by the equation:
$$Y^2+(X^2-2)(X^2-1)(X^2+1)=0,$$
endowed with the action of $\Z/2\Z$ given by $(X,Y) \mapsto (-X,Y)$. The curve $C$ consists of two connected components homeomorphic to $S^1$, that are exchanged under the action of $\Z/2\Z$. The algebraic quotient is the real algebraic irreducible plane curve $\widetilde C$:
$$T^2+(U-2)(U-1)(U+1)=0.$$
It contains the topological quotient of $C$ as a connected component homeomorphic to $S^1$, and an extra connected component coming from the conjugated complex points of $C$.
\end{rmk}
%%%%%%%%%%%%%%%%%%%%%%%%%%%%%%%%%%%%%%%%%%%%%%%%%%%%%%%%%%%%%%%%%%%%%%%%%%
%%%%%%%%%%%%%%%%%%%%%%%%%%%%%%%%%%%%%%%%%%%%%%%%%%%%%%%%%%%%%%%%%%%%%%%%%%
\section{The case $G=\Z/2\Z$}\label{G=2}

The case $G=\Z/2\Z$ is of great interest with a view to applications in section \ref{bN}. We pay particular attention to understanding the negative part of the \veps in accordance with the non-virtual case (cf. section \ref{negative}). We concentrate also on necessary and sufficient conditions for the quotient of an arc-symmetric set to remain arc-symmetric (cf. \cite{hui} also).

\subsection{Quotient of arc-symmetric sets by $\Z/2\Z$-action}

We focus on the quotient of an arc-symmetric set by the action of $G=\Z/2\Z$. In particular, we give necessary and sufficient conditions on the action in order that the quotient remains arc-symmetric.

In general, the quotient space of a real algebraic set is no longer an algebraic set, but only a semi-algebraic set. More precisely, denote by $X$ a real algebraic set and by $\sigma : X \dans X$ an involution. The quotient of the complexification $X_{\C}$ of $X$ by the complexified action $\sigma_{\C}$ of $\sigma$ is a complex algebraic variety $Y$ endowed with an antiholomorphic involution that commutes with $\sigma$. The image of $X$ in $Y$ is a semi-algebraic subset of the real algebraic set $Y(\R)$ of fixed points of $Y$ under the antiholomorphic involution. By restriction, the same situation holds concerning arc-symmetric sets with an algebraic action. We are interested in the case when this quotient remains arc-symmetric.

By the quotient $B$ of an arc-symmetric set $A$, we mean a surjective map $p:A \dans B$ where $p$ is the restriction of a regular map from the Zariski closure $\overline A^Z$ of $A$, with the property that $p(x)=p(y)$ if and only if $y\in \{ x, \sigma (x)\}$. In proposition \ref{}, we give necessary and sufficient conditions on the action for the semi-algebraic set $B$ to be arc-symmetric.

In case $X$ is nonsingular (not necessarily compact) but the quotient space is arc-symmetric, then the action on an equivariant nonsingular compactification must be free, unless it is trivial on a whole component. Consider first the compact case.

\begin{prop}\label{quot-free} Let $X$ be a compact nonsingular arc-symmetric set endowed with a $G$-action. Then the quotient space $X/G$ is arc-symmetric if and only if the action of $G$ on $X$ is trivial on each irreducible component of $X$ containing a fixed point. 
\end{prop}

\begin{proof} One can assume $X$ to be algebraic without loss of generality. Denote by $p$ the quotient morphism $p:X\dans X/G$. Denote by $p_{\C}:X_{\C}\dans X_{\C}/G$ its complexification. Then $X/G$ is a semi-algebraic subset of the real points of the complex algebraic variety $X_{\C}/G$.

Assume that $x \in X$ is a fixed point under the action of $G$. By lemma \ref{inter-lem} below, there exists an injective analytic arc $\delta :]-1,1[\dans X$ with origin at $x$ so that $\delta (-t)=\sigma \delta (t)$. Then the analytic arc $p(\delta)$ satisfies $p(\delta)(-t)=p(\delta)(t)$, and therefore there exists an analytic arc $\gamma$ so that $p(\delta)(t)=\gamma(t^2)$. Remark that $\gamma(]0,1[)$ is included in $X/G$ whereas $\gamma(]-1,0[)$ does not even meet $X/G$ since $\gamma(-t^2)=p(\delta)(it)\in (X_{\C}/G)\setminus (X/G)$. This prove that the quotient space $X/G$ can not be arc-symmetric.
\end{proof}

Let us prove now the intermediate result.

\begin{lemma}\label{inter-lem} Let $X$ be a nonsingular arc-symmetric set endowed with a nontrivial involution $\sigma$ that admits a fixed point $x\in X$. Then there exists an injective analytic arc $\delta :]-1,1[\dans X$ with origin at $x$ so that $\delta (-t)=\sigma \delta (t)$.
\end{lemma}

\begin{proof} The action of $\sigma$ on $X$ is locally linearisable around the fixed point $x$ by \cite{cartan}: there exist an open neighbourhood $U$ of $x$ in $X$ stable under $\sigma$, an open neighbourhood $V$ of $0$ in $\R^{\dim X}$, and an analytic isomorphism $f:U \dans V$ so that $f(x)=0$ and the action of $\sigma$ on $U$ is transported via $f$ in a linear action on $V$. Now, one can choose a suitable analytic arc $\gamma$ on $V$ through $0$ that lies in the eigenspace associated to the eigenvalue $-1$ to obtain the expected arc $\delta=f^{-1}(\gamma)$ on $U$.
\end{proof}

If the arc-symmetric set is no longer nonsingular, then the equivalence of proposition \ref{quot-free} fails as illustrated by the last example in \ref{exe}.

However for a nonsingular arc-symmetric set $X$, we obtain a necessary and sufficient condition for the existence of an arc-symmetric quotient of $X$ in terms of the action on a nonsingular equivariant compactification of $X$.

\begin{prop}\label{eq-lisse} Let $X$ be a nonsingular arc-symmetric set endowed with a $G$-action. Let $\widetilde X$ be a nonsingular equivariant compactification of $X$ with $\widetilde X\setminus X$ compact. Then the quotient space $X/G$ is arc-symmetric if and only if the action of $G$ on $\widetilde X$ is trivial on each irreducible component of $\widetilde X$ containing a fixed point.
\end{prop}

\begin{proof} If $X/G$ is arc-symmetric, then the proof of proposition \ref{quot-free} implies that the action of $G$ on $\widetilde X$ is free, or trivial on some components. Conversely the semi-algebraic quotient of $X$ by $G$ can be decomposed as the difference between $\widetilde X/G$ and $(\widetilde X \setminus X)/G$. The former quotient is arc-symmetric by proposition \ref{quot-free}, the latter by proposition \ref{action-free}. Therefore $X/G$ is also arc-symmetric because arc-symmetric sets are stable under boolean operations.
\end{proof}

%%%%%%%%%%%%%%%%%%%%%%%%%%%%%%%%%%%%%%%%%%%%%%%%%%%%%%%%%%%%%%%%%%%%%%%%%%%
\subsection{$\Z/2\Z$-equivariant virtual Betti numbers}

If we can not relate the equivariant virtual Poincar\'e series of $X$ endowed with a free $G$-action to the virtual Poincar\'e polynomial of the quotient in general, we can determine the negative part of the equivariant virtual Poincar\'e series in the case $G=\mathbb Z/2$.

\begin{lemma}\label{negative-serie} Let $X$ be an arc-symmetric set endowed with a free $G$-action. Then $\beta_n^G(X)=0$ for $n <0$.
\end{lemma}

Before beginning with the proof, let us recall that in case $X$ is compact and nonsingular then so is the fixed point set $X^G$.

\begin{proof} We make an induction on the dimension. Note that one can assume $X$ to be nonsingular by additivity of the equivariant virtual Poincar\'e series. Now, take a nonsingular equivariant compactification $X \subset \overline X$ so that the dimension of $\overline X \setminus X$ is strictly less than that of $X$. The compactification may have fixed points, nevertheless the fixed points set $ \overline X ^G$ is nonsingular. Now
$$\beta_n^G(X)=\beta_n^G(\overline X)-\beta_n^G(\overline X \setminus X)$$
by additivity and $$\beta_n^G(\overline X)=b_n^G(\overline X)=\sum_{i\geq 0}b_i(\overline X^G)$$
since $\overline X$ is compact and nonsingular (cf. section \ref{negative}). However, by the induction assumption
$$\beta_n^G(\overline X \setminus X)=\sum_{i\geq 0}\beta_i((\overline X \setminus X)^G).$$
Since the fixed points sets of $\overline X$ and $\overline X \setminus X$ coincide, the right hand side reduces to
$$ \sum_{i\geq 0}\beta_i(\overline X^G).$$
We complete the proof by recalling that $\overline X^G$ is compact and nonsingular therefore $\beta_i(\overline X^G)=b_i(\overline X^G)$ by definition of the virtual Poincar\'e polynomial.
\end{proof}

Recall that the negative equivariant homology groups can be computed in terms of the homology groups of the fixed points (section \ref{negative}) when the group acting is $G=\Z/2\Z$. Similarly we obtain:

\begin{prop}\label{neg} Let $X$ be an arc-symmetric set with an action of $G=\Z/2\Z$. Take $n<0$. Then $\beta_n^G(X)$ is equal to the virtual Poincar\'e polynomial of $X^G$ evaluated at $1$:
$$\beta_n^G(X)=\sum_{i \geq 0}\beta _i (X^G).$$
\end{prop}

\begin{proof} Note that
$$\beta^G(X)=\beta^G(X^G)+\beta^G(X \setminus X^G)=\beta(X^G)\beta^G(\textbf{pt})+\beta^G(X \setminus X^G),$$
therefore the proof is a direct consequence of lemma \ref{negative-serie} and of the computation of $\beta^G(\textbf{pt})$ in remark \ref{point}.
\end{proof}

\begin{ex}\label{exe} Consider the curve $C$ given by the equation $Y^2=X^2-X^4$ in $\mathbb R^2$ endowed with different actions of $\Z/2\Z$. Note that the resolution of the singularities of $C$ obtained by blowing-up the unique singular point $p=(0,0)$ is homeomorphic to $S^1$. The inverse image of this singular point $\{p\}$ consists of two points, denoted by $\{p_1,p_2\}$, either fixed or inverted by the action. Therefore by additivity one can express the \veps of $C$ by:
$$\beta^G(C)=\beta^G(S^1)-\beta^G(\{p_1,p_2\})+\beta^G(\{p\}).$$
\begin{itemize}
\item $(X,Y) \mapsto (-X,-Y)$. The fixed points of the resolution are exactly $p_1$ and $p_2$, thus:
$$\beta^G(C)=u+1+\frac{1}{u-1}.$$
\item $(X,Y) \mapsto (X,-Y)$. The action inverts the points $p_1$ and $p_2$, but is not free on the resolution, so:
$$\beta^G(C)=u+2+\frac{3}{u-1}.$$
\item $(X,Y) \mapsto (-X,Y)$. The action inverts the points $p_1$ and $p_2$, and is free on the resolution, therefore by example \ref{ex-beta}.\ref{point-ech}:
$$\beta^G(C)=u+1+\frac{1}{u-1}.$$
Note that in this example the algebraic quotient of $C$ does exist and is given by the circle defined by $Y^2=X-X^2$ in $\mathbb R^2$. Moreover $\beta^G(C) \neq \beta(C/G)$ and therefore proposition \ref{action-free} does not hold in case the action is no longer free.
\end{itemize}
\end{ex}

\begin{prop}\label{cle} Let $X$ be a nonsingular arc-symmetric set endowed with an action of $G=\Z/2\Z$ so that the quotient space is arc-symmetric. Assume moreover that the action is not trivial on any irreducible component of $X$. Then $$\beta^G(X)=\beta(X/G).$$
\end{prop}

\begin{proof} We prove the result by induction on the dimension. Let $\widetilde X$ be a nonsingular equivariant compactification of $X$. Then the $G$-action on $\widetilde X$ is free and the quotient space $\widetilde X/G$ is arc-symmetric by proposition \ref{eq-lisse}. By additivity of $\beta^G$:
$$\beta^G(X)=\beta^G(\widetilde X)- \beta^G(\widetilde X \setminus X).$$
Moreover $\beta^G(\widetilde X)=\beta(\widetilde X/G)$ by proposition \ref{action-free} whereas $\beta^G(\widetilde X\setminus X)=\beta((\widetilde X\setminus X)/G)$ by the induction assumption. Finally
$$\beta^G(X)=\beta(\widetilde X/G)-\beta((\widetilde X\setminus X)/G)=\beta(X/G)$$
by additivity of $\beta$.
\end{proof}
%%%%%%%%%%%%%%%%%%%%%%%%%%%%%%%%%%%%%%%%%%%%%%%%%%%%%%%%%%%%%%%%%%%%%%%%%%
%%%%%%%%%%%%%%%%%%%%%%%%%%%%%%%%%%%%%%%%%%%%%%%%%%%%%%%%%%%%%%%%%%%%%%%%%%
%%%%%%%%%%%%%%%%%%%%%%%%%%%%%%%%%%%%%%%%%%%%%%%%%%%%%%%%%%%%%%%%%%%%%%%%%%
\section{Application to blow-Nash equivalence}\label{bN}
The study of real germs is difficult in general, notably in the choice of a good equivalence relation between these germs. We focus here on the blow-Nash equivalence between Nash function germs. It is an analog of the blow-analytic equivalence defined by T.C. Kuo \cite{Kuo} for Nash function germs. Let us mention in particular that there are no moduli for families with isolated singularities \cite{fichou}. We know also invariants that enable us to begin some classification results \cite{fichou1}.

As an application of blow-Nash equivalence, we introduce equivariant zeta functions of a Nash function germ defined via the equivariant virtual Poincar\'e series. They generalise the zeta functions introduced in \cite{fichou} with the virtual Poincar\'e polynomial (see also \cite{KP}). We prove in particular their rationality (proposition \ref{DLmono}), and their invariance under the blow-Nash equivalence between Nash function germs. We achieve moreover a better understanding of the zeta functions of constant sign Nash germs. 

\subsection{Equivariant Zeta functions}\label{zeta}
We define the equivariant zeta functions for a germ of Nash
functions in the spirit of \cite{DL1}. A Nash function germ $f:(\mathbb R ^d,0) \longrightarrow (\mathbb R,0)$ is a real analytic germ with semi-algebraic graph.

Denote by $\mathcal L$ the space of arcs at the
origin $0 \in \mathbb R ^d$, defined by:
$$\mathcal L=\mathcal L(\mathbb R ^d,0)= \{\gamma : (\mathbb R,0) \longrightarrow (\mathbb R ^d,0)
:\gamma \textrm{ formal}\},$$
and by $\mathcal L_n$ the space of truncated arcs at the order $n+1$:
$$\mathcal L_n=\mathcal L_n(\mathbb R ^d,0)= \{\gamma \in \mathcal L:
\gamma (t)=a_1t+a_2t^2+ \cdots a_nt^n,~a_i \in \mathbb R ^d\},$$
for $n\geq 0$ an integer.

Consider a Nash function germ $f:(\mathbb R ^d,0) \longrightarrow (\mathbb R,0)$. We define the equivariant zeta functions of $f$ to be the formal power series with coefficients in $\mathbb Z[[u^{-1}]][u]$:
$$Z_{f,+}^{G}(T)= \sum _{n \geq 1}{\beta ^G(A_n^{+})u^{-nd}T^n}
\textrm{~~and~~ }Z_{f,-}^{G}(T)= \sum _{n \geq 1}{\beta^G (A_n^{-})u^{-nd}T^n},$$
where 
$$A_n^{+} =\{\gamma \in  \mathcal L_n: f\circ \gamma
(t)=+t^n+\cdots \}\textrm{~~and~~ } A_n^{-} =\{\gamma \in  \mathcal L_n: f\circ \gamma
(t)=-t^n+\cdots \} $$
are endowed with the action of $G=\Z/2\Z$ defined by
\begin{itemize}
\item $t \mapsto t$ if $n$ is odd,
\item $t \mapsto -t$ if $n$ is even.
\end{itemize}

Note that the spaces of truncated arcs $A_n^{+}$ and $A_n^{-}$ for $n \geq 1$ are Zariski constructible subsets of
$\mathbb R^{nd}$ for $n \geq 1$. In particular they belong to $\mathcal{AS}$ and their equivariant virtual Poincar\'e series are well-defined.

\begin{rmk}\label{remark} One recovers the zeta functions with sign introduced in \cite{fichou} by considering the action of the trivial group on the arc spaces $A_n^{+}$ and $A_n^{-}$.
\end{rmk}

 Let $f:(\mathbb R^d,0) \longrightarrow
  (\mathbb R,0)$ be a Nash function germ.
Let $h:\big(M,h^{-1}(0)\big)  \longrightarrow (\mathbb R ^d,0)$ be a proper surjective Nash map such that $f \circ h$ and the
Jacobian determinant $\jac h$ have simultaneous normal crossings. Assume moreover that $h$ is an isomorphism over the
complement of the zero locus of $f$.

We need some notation in order to express the equivariant zeta functions of $f$ in terms of the Nash modification $h$ of $f$. 
Let $(f \circ h)^{-1}(0)= \cup_{j \in J}E_j$ be the decomposition
into irreducible components of $(f \circ h)^{-1}(0)$, and assume
that $ h^{-1}(0)=\cup_{k \in K}E_k$ for some $K \subset J$.

Put $N_i=\mult _{E_i}f \circ h$ and $\nu _i=1+\mult _{E_i} \jac
h$, and for $I \subset J$ denote by $E_I^0$ the set $(\cap _{i
\in I} E_i) \setminus (\cup _{j \in J \setminus I}E_j)$.

One defines a covering $\widetilde
{E_I^{0,\pm}}$ of $E_I^0$, where $\pm$ designs either $+$ or $-$, in
the following way. For $U$ an affine open subset of $M$
such that $f \circ h=u \prod_{i\in I}y_i^{N_i}$ on $U$, where $u$
is a Nash unit, we put 
$$R_{U}^{\pm}=\{ (x,s) \in (E_I^0 \cap U) \times \mathbb
R; s^{m_I}=\pm \frac{1}{u(x)}\}$$ where $m_I=gcd(N_i)$. Then the sets $R_{U}^{\pm}$ glue
together along the $E_I^0 \cap U$ to define $\widetilde {E_I^{0,\pm}}$.

The covering $\widetilde {E_I^{0,\pm}}$ is endowed with the trivial $G$-action in the case $m_I$ is odd, and by the $G$-action induced by the multiplication by $-1$ on the variable $s$ in the case $m_I$ is even.

Now we can state the rationality of the equivariant zeta functions of a Nash function germ. The formula below is called the Denef and Loeser formula.

\begin{prop}\label{DLmono}  The equivariant zeta functions of a Nash function germ are rational. More precisely, with the assumptions and notations above, the equivariant zeta functions equals:

$$Z^G_{f,{\pm}}(T)=\sum_{I\neq \emptyset} (u-1)^{|I|-1}\beta ^G\big(\widetilde{E_I^{0,\pm}} \cap
h^{-1}(0)\big) \prod_{i \in I}\frac{u^{-\nu_i}T^{N_i}}{1-u^{-\nu_i}T^{N_i}}.$$
\end{prop}

\begin{proof} Thanks to proposition \ref{mult} the proof of the Denef and Loeser formula in the non-equivariant case (\cite{fichou2}) carries over to the equivariant case. In particular, the $G$-action introduced on the covering $\widetilde {E_I^{0,\pm}}$ of $E_I^0$ is actually induced by the $G$-action on the arc spaces $A_n^{\pm}$. 
\end{proof}

\begin{ex}\label{ex-pos} Let $f:\mathbb R^2 \longrightarrow
    \mathbb R$ be defined by $f(x,y)=x^2+y^4$. One can resolve the
    singularities of $f$ by two successive blowings-up. The exceptional divisor $E$ has 
    two irreducible components $E_1$ and $E_2$ with $N_1=2,~ \nu_1=2,
   ~ N_2=4, ~\nu _2=3$. Moreover $\widetilde {E_{\{1\}}^{0,+}}$  and
    $\widetilde {E_{\{2\}}^{0,+}}$ are homeomorphic to a circle minus
    two points and the $G$-action has fixed points on the closure. Therefore
$$Z_{f,+}^G(T)=(u-1)\frac{u^{-2}T^2}{1-u^{-2}T^2}\frac{u^{-3}T^4}{1-u^{-3}T^4}+u
\frac{u^{-2}T^2}{1-u^{-2}T^2}+u \frac{u^{-3}T^4}{1-u^{-3}T^4}.$$
\end{ex}
%%%%%%%%%%%%%%%%%%%%%%%%%%%%%%%%%%%%%%%%%%%%%%%%%%%%%%%%%%%%%%%%%%%%%%%%%%
\subsection{Blow-Nash equivalence}
Let us state the definition of the blow-Nash equivalence between
Nash function germs. It consists of an adaptation of the blow-analytic equivalence defined by T.-C. Kuo (\cite{Kuo}) to the Nash framework.
\begin{defi}\label{defbN}\begin{flushleft}\end{flushleft}
\begin{enumerate}
\item Let $f:(\mathbb R ^d,0) \longrightarrow (\mathbb R,0)$ be a Nash
  function germ. A Nash modification of $f$ is a proper
  surjective Nash map $\pi:\big( M,\pi^{-1}(0) \big)
  \longrightarrow (\mathbb R^d,0)$ whose complexification $\pi^*$ is an isomorphism everywhere except on some thin subset of $M^*$, and 
such that $f \circ \pi$ has only normal crossings.
\item \label{point2} Two given germs
of Nash functions $f,g:(\mathbb R ^d,0) \longrightarrow (\mathbb R,0)$ are said to be blow-Nash equivalent if there
exist two Nash modifications $$h_f~:~\big(M_f,h_f^{-1}(0)\big)  \longrightarrow
(\mathbb R ^d,0) \textrm{ and } h_g:\big(M_g,h_g^{-1}(0)\big)  \longrightarrow
(\mathbb R ^d,0),$$ and a Nash isomorphism $\Phi$
between semi-algebraic neighbourhoods $\big(M_f,h_f^{-1}(0)\big)$ and $\big(M_g,h_g^{-1}(0)\big)$ that preserves the
  multiplicities of the Jacobian determinant along the exceptional divisors
of the Nash modifications
  $h_f,~h_g$, and that
induces a homeomorphism $\phi :(\mathbb R ^d,0) \longrightarrow
(\mathbb R^d,0)$ such that the diagram
$$\xymatrix{\big(M_f,h_f^{-1}(0)\big) \ar[rr]^{\Phi} \ar[d]_{h_f}& &\big(M_g,h_g^{-1}(0)\big) \ar[d]^{h_g}\\
            (\mathbb R ^d,0)       \ar[rr]^{\phi} \ar[dr]_f        & &           (\mathbb R ^d,0) \ar[dl]^g \\
         &   (\mathbb R,0) & } $$
is commutative.
\end{enumerate}
\end{defi}

\begin{prop} Let $f,g:(\mathbb R ^d,0) \longrightarrow (\mathbb R,0)$ be blow-Nash equivalent Nash function germs. Then the equivariant zeta functions of $f$ and $g$ coincide.
\end{prop}

The proof is similar to that of the non-equivariant case \cite{fichou}, using propositions \ref{nash-iso} and \ref{DLmono}.
%%%%%%%%%%%%%%%%%%%%%%%%%%%%%%%%%%%%%%%%%%%%%%%%%%%%%%%%%%%%%%%%%%%%%%%%%%
\subsection{Constant sign Nash function germs}
%%%%%%%%%%%%%%%%%%%%%%%%%%%%%%%%%%%%%%%%%%%%%%%%%%%%%%%%%%%%%%%%%%%%%%%%%%
In \cite{fichou}, we introduced zeta functions with signs, which correspond to the equivariant zeta functions considered with the trivial $G$-action as noticed in remark \ref{remark}. We introduced also a ``naive'' zeta function, defined by:
$$Z_{f}(T)= \sum _{n \geq 1}{\beta ^G(A_n)u^{-nd}T^n}$$
where 
$$A_n =\{\gamma \in  \mathcal L_n: f\circ \gamma
(t)=ct^n+\cdots,~~c\neq 0 \}.$$
It was not clear what relations hold between these zeta functions, even in the case of constant sign Nash germs. For example, compare the positive zeta function of example \ref{ex-pos} with the ``naive'' zeta function:
$$Z_f(T)=(u-1)^2\frac{u^{-2}T^2}{1-u^{-2}T^2}\frac{u^{-3}T^4}{1-u^{-3}T^4}+(u-1)u
\frac{u^{-2}T^2}{1-u^{-2}T^2}+(u-1)u \frac{u^{-3}T^4}{1-u^{-3}T^4}$$
and the positive non-equivariant zeta functions of \cite{fichou}:
$$Z_f^+(T)=2(u-1)\frac{u^{-2}T^2}{1-u^{-2}T^2}\frac{u^{-3}T^4}{1-u^{-3}T^4}+(u-1)
\frac{u^{-2}T^2}{1-u^{-2}T^2}+(u-1) \frac{u^{-3}T^4}{1-u^{-3}T^4}.$$

The equivariant zeta functions solve this issue.

\begin{prop} Let $f:(\mathbb R^d,0) \longrightarrow (\mathbb R,0)$ be a nonnegative Nash function germ. Then $$Z_{f}(T)=(u-1)Z^G_{f,+}(T).$$
\end{prop}

\begin{proof} The expression for the naive zeta functions via its Denef and Loeser formula is \cite{fichou}:
$$Z_f(T)=\sum_{I\neq \emptyset} (u-1)^{|I|}\beta\big(E_I^0 \cap
h^{-1}(0)\big) \prod_{i \in I}\frac{u^{-\nu_i}T^{N_i}}{1-u^{-\nu_i}T^{N_i}} $$ 
with the notations of proposition \ref{DLmono}.
Therefore, comparing to the Denef and Loeser formula for $Z^G_{f,+}(T)$ of proposition \ref{DLmono}, it suffices to prove that:
$$\beta^G(R^+_U)=\beta(E_I^0 \cap U)$$ 
The result follows from proposition \ref{cle} since $E_I^0 \cap U$ is the quotient of $R^+_U$ under the $\Z/2\Z$-action.
\end{proof}
%%%%%%%%%%%%%%%%%%%%%%%%%%%%%%%%%%%%%%%%%%%%%%%%%%%%%%%%%%%%%%%%%%%%%%%%%%
%%%%%%%%%%%%%%%%%%%%%%%%%%%%%%%%%%%%%%%%%%%%%%%%%%%%%%%%%%%%%%%%%%%%%%%%%%
%%%%%%%%%%%%%%%%%%%%%%%%%%%%%%%%%%%%%%%%%%%%%%%%%%%%%%%%%%%%%%%%%%%%%%%%%%

\end{document}